# ON THE GENERALIZED RIEMANN HYPOTHESIS II

D. Ghisa

York University, Glendon College, Toronto, Canada
dghisa@yorku.ca

*Some of my previous publications were incomplete in the sense that non trivial zeros belonging to a particular type of fundamental domain have been inadvertenly ignored. Due to this fact, I was brought to believe that computations done by some authors in order to show counterexamples to RH were affected of approximation errors. In this paper I illustrate graphically the correctness of those computations and I fill the gaps in my publications.*

## 1. INTRODUCTION

In his 1934 Ph.D. thesis [21] A. Utzinger introduces a new method of study of the Riemann Zeta function by using the pre-image by $\zeta$ of the real axis (*reelen Züge*).The method looked promising since in the same year by using the respective method A. Speiser (see [19]) succeeded to prove an equivalent form of the Riemann Hypothesis (RH). Speiser's result is mentioned in most of the works in the field although, as Arias-de-Reyna [5] noticed, "nobody reproduces his theorems". The reason seems to be a certain reluctance regarding the argumentation which uses geometry. However, using conformal mapping techniques, which are geometric in their nature, does not seem to produce any adverse reaction. This is why, starting with the 2009 paper [2] we have tried to use only arguments strictly related to the conformal mapping.

The concept of fundamental domain coming from the automorphic functions theory and to which Ahlfors [1] gave an independent meaning helped us to fulfill this task and to advance the research in the field of Dirichlet series. The Utzinger and Speiser *reellen Züge* were used only to put into evidence in [3], [15], [12], [13] the boundaries of fundamental domains of the class of functions we were studying. This class of functions includes the Selberg class, for which the most general RH has been formulated. An important step forward was the discovery of the strips $S_k$ (described in section 2) allowing the use of a divide and conquer algorithm. It was now enough to deal with just an arbitrary strip $S_k$ in order to draw conclusions valid for the whole complex plane. Yet, the configurations inside a strip $S_k$ needed a geometric way of thinking and they eluded us with a few occasions.

The regularity of the curves $\Gamma_{k,j}$ and $\Upsilon_{k,j}$ (see [12] and [13]) expressed by two elementary topological properties, the *color alternating rule* and the *color matching rule* is subject also to a variety of configurations which were impossible to imagine at the beginning. *Embracing curves* came as a surprise for high values of $t$, only to notice later that it is a common phenomenon and a curve $\Gamma_{k,j}$, $j \neq 0$ can embrace several such curves, but never a curve $\Gamma_{k,0}$. Indeed, every $\Gamma_{k,0}$ extends for $\sigma$ in the whole range of real numbers, while the curves $\Gamma_{k,j}$, $j \neq 0$ remain bounded to the right.

However, as seen in the case of Davenport and Heilbronn function, $\Gamma_{k,0}$ and a component of the pre-image of a close ray can behave as an embracing curve. The embracing effect on the corresponding fundamental domains is that they become bounded to the right. It was impossible to realize that a similar effect can produce some curves $\Gamma_{k,0}$. Only by an appropriate zooming on the zeros of Davenport and Heilbronn function (see [20], [6], [8], [5]) we succeeded to illustrate this fact and to notice that the exceptions to RH happen exactly with the zeros located on some curves $\Gamma_{k,0}$.

The study of conformal mappings of the fundamental domains appears to be of a crucial importance when dealing with the distribution of non trivial zeros of the functions, since every fundamental domain contains a unique such zero and the location of the respective zero is conditioned by the corresponding conformal mapping.

The strips $S_k$ encapsulate deep properties of the classes of functions we are studying. Everyone of these strips contains a unique unbounded component of the pre-image of the unit disk, as well as a unique component $\Gamma_{k,0}$ of the pre-image of the real axis which is mapped bijectively by the function onto the interval $(-\infty, 1)$ of the real axis.

For the Riemann Zeta function the width of every strip $S_k$ in the neighborhood of the imaginary axis is roughly of 9 units, as we could check for values of the ordinate $t$ up to $10^9$.

The number of zeros in $S_k$ increases with $|t|$ following approximately a logarithmic law.

In Fig. 1 a zoom is shown on the first four couple of symmetric zeros with respect to the critical line of the Davenport and Heilbronn function from the list presented in [6].

It gives us a pretty clear idea of how the conformal mappings of the corresponding fundamental domains look like in the neighborhood of these points.

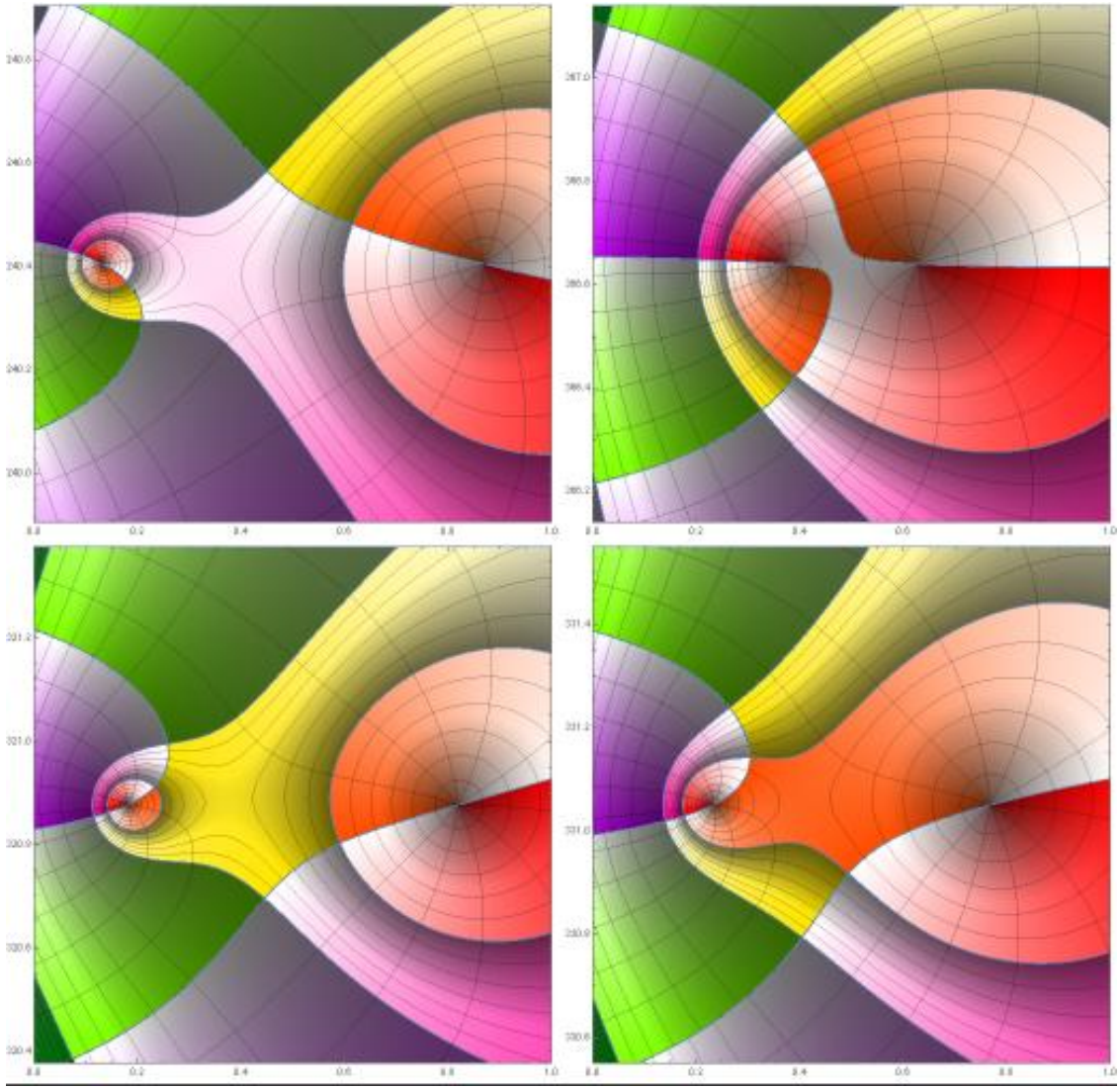

Fig.1 Off critical line non trivial zeros of Davenport and Heilbronn function

Details of these conformal mappings are illustrated in Fig.2. Here we take a closer look to the first and the last picture. They portray the conformal mappings in the neighborhood of symmetric zeros with respect to the critical line located on $\Gamma_{k,0}$ and $\Gamma_{k,1}$ and respectively $\Gamma_{k,0}$ and $\Gamma_{k,-1}$. Those zeros have the abscissa $0.86953\ldots$ and $0.13046\ldots$, respectively $0.76822\ldots$ and $0.23177\ldots$

By a rough computation of the corresponding zeros of the derivative we have found the abscissa $0.31\ldots$ and respectively $0.39\ldots$

There are two important facts to be noticed about these last points: the first one is that they are located on the left of the critical line and the second is that each one has the same imaginary part as the two corresponding zeros of the function.

Are these facts a simple coincidence, or do they represent necessary conditions? We will answer this question in section 3.

We were wrong assuming that the deviation from the critical line of Spira's and Balazario's points is due to errors of approximation (see [10], [14] ) and they are not true non trivial zeros. This assumption was supported by the fact that in [10] we took inadvertently a wrong Dirichlet character in building what we considered to be the Davenport and Heilbronn function.

However, we are reluctant to say that those zeros represent *counterexamples* to RH, since they are zeros of a function which is not a Dirichlet L-function neither is it a member of the Selberg's class of functions.

Now we need to prove that the unwanted behavior of $\Gamma_{k,0}$ we have noticed in this case cannot happen in the case of Dirichlet L-functions and to complete in this way the proof of the Theorem 9 from [12], the description given in [7] for the location of the zeros of the derivative of Dirichlet L-functions, as well as the Section 9.3 from [13] .

On the other hand a more general situation has been considered in [14] and we have to deal here with generalized Dirichlet series.

The monograph [16] was probably the first dedicated entirely to the general theory of Dirichlet series. An elegant exposure of this topic, as well as interesting applications to Bohr functions can be found in [4]. Yet, the motivation for such an approach comes from a lot of other recent developments far beyond the scope of this paper.

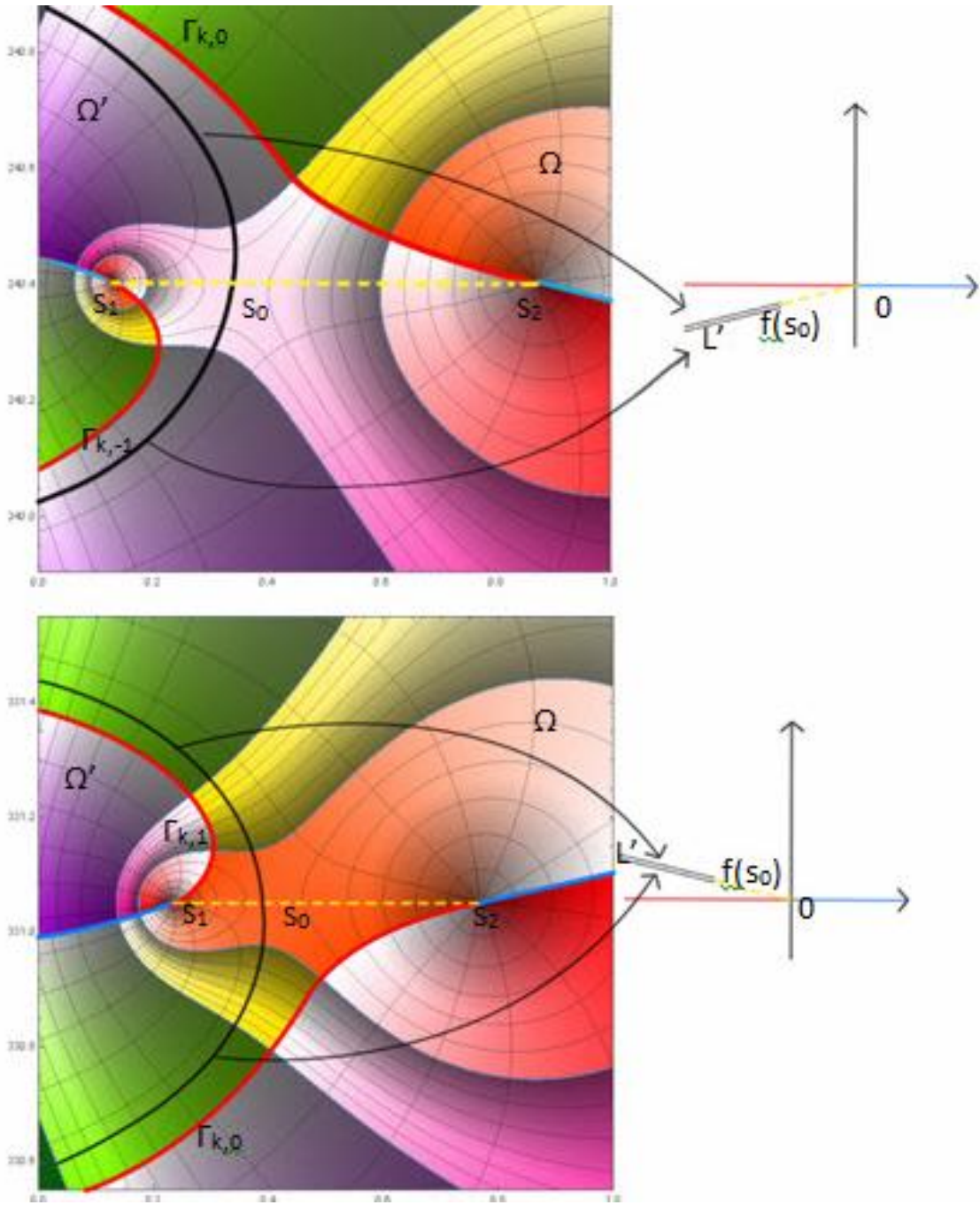

Fig. 2 Conformal mapping of the fundamental domains.

## 2. GENERAL DIRICHLET SERIES DEFINED BY TOTALLY MULTIPLICATIVE FUNCTIONS

What makes the Davenport and Heilbronn function special in contrast with Dirichlet L-functions is the fact that although it satisfies a Riemann type of functional equation, the corresponding Dirichlet series is not an Euler product. However, it has been obtained as a linear combination of series which were Euler products. We have shown in [11] that this function is not so special after all, since infinitely many functions of this type can be built using appropriate Dirichlet L-series. The question remains if all of them violate RH.

Since we do not know the answer to this question, for the purpose of dealing with functions satisfying RH, we need to limit our attention to Dirichlet series which can be expressed as Euler products, although as shown in [9], this is not a necessary condition for the validity of RH. We studied in [14] general Dirichlet series (see [16] and [4], Chapter 8) of the form:

$$\zeta_{A,\Lambda}(s) = \Sigma_{n=1}^{\infty} a_n e^{-\lambda_n s} \qquad (1)$$

where $A = \{a_n\}_{n \in \mathbb{N}}$ is an arbitrary infinite sequence of complex numbers and $\Lambda = \{0 = \lambda_1 \leq \lambda_2 \leq \ldots\}$ is a non decreasing sequence of non negative numbers such that $\lim_{n \to \infty} \lambda_n = +\infty$. We suppose that the abscissa of convergence $\sigma_c$ of the series (1):

$$\sigma_c = \limsup\nolimits_{n \to \infty} \ln|\Sigma_{k=1}^{n} a_k|^{1/\lambda_n} \qquad (2)$$

is finite and that this series can be continued analytically to the whole complex plane except for a possible pole at $s = 1$. We keep the notation $\zeta_{A,\Lambda}(s)$ for the respective function. Since $\lim_{\sigma \to +\infty} \zeta_{A,\Lambda}(\sigma + it) = 1$, uniformly with respect to $t$ (see [14]) there is $\sigma_0 \geq \sigma_c$ such that $\zeta_{A,\Lambda}(s) \neq 0$ for $\operatorname{Re} s > \sigma_0$.

Suppose that $\zeta_{A,\Lambda}(s)$ satisfies a Riemann-type of functional equation:

$$\zeta_{A,\Lambda}(s) = M(s)\overline{\zeta}_{A,\Lambda}(1-s) \qquad (3)$$

where $\overline{\zeta}_{A,\Lambda}(s) = \overline{\zeta_{A,\Lambda}(\bar{s})})$ and $M(s)$ is an analytic function in $\mathbb{C}$. Then the zeros of $M(s)$ will be called *trivial zeros* of $\zeta_{A,\Lambda}(s)$ and those zeros of $\zeta_{A,\Lambda}(s)$ which are not zeros of $M(s)$ will be called *non trivial zeros*. This convention is not quite conformal with the traditional one (see [17], page 333), yet in dealing with general Dirichlet series we feel that it is more convenient to stick with it. We notice that, due to the functional equation (3), if $\zeta_{A,\Lambda}(\sigma_0 + it_0) = 0$ and $M(\sigma_0 + it_0) \neq 0$, then necessarily $\zeta_{A,\Lambda}(1 - \sigma_0 + it_0) = 0$. The RH for a given class of functions $\zeta_{A,\Lambda}(s)$ says that this happens in that class if and only if $\sigma_0 + it_0 = 1 - \sigma_0 + it_0$, i.e. $\sigma_0 = 1/2$.

The function $n \to a_n$, $a_n \in A$ is called totally multiplicative if $a_1 = 1$ and for any $j, k \in \mathbb{N}$ we have $a_{jk} = a_j a_k$. If $n \to a_n$ is totally multiplicative, then for any $n$ with the prime decomposition $n = p_1^{\alpha_1} p_2^{\alpha_2} \ldots p_k^{\alpha_k}$ we have $a_n = a_{p_1}^{\alpha_1} a_{p_2}^{\alpha_2} \ldots a_{p_k}^{\alpha_k}$. From now on we will suppose that $n \to a_n$ is a totally multiplicative function and also that the function $n \to \lambda_n$ has the property that

$$n = p_1^{\alpha_1} p_2^{\alpha_2} \ldots p_k^{\alpha_k} \Rightarrow \lambda_n = \alpha_1 \lambda_{p_1} + \alpha_2 \lambda_{p_2} + \ldots + \alpha_k \lambda_k. \quad (4)$$

We notice that this property holds when $\lambda_n = \ln n$, i.e. when (1) is an ordinary Dirichlet series and also that $a_n = \chi(n)$ is a totally multiplicative function when $\chi$ is a Dirichlet character of an arbitrary modulus $q$, i.e. the Dirichlet L-functions belong to the class we are dealing with.

**Theorem 1.** *If $n \to a_n$ is a totally multiplicative function and $n \to \lambda_n$ satisfies the condition (4), then for $\operatorname{Re} s > \sigma_c$*

$$\zeta_{A,\Lambda}(s) = \prod_p (1 - a_p e^{-\lambda_p s})^{-1} \quad (5)$$

*where $p$ runs through the prime numbers.*

*Proof:* Let us write (1) under the form: $\zeta_{A,\Lambda}(s) = 1 + a_2 e^{-\lambda_2 s} + a_3 e^{-\lambda_3 s} + \ldots$

Then, for $\operatorname{Re} s \geq \sigma_0 > \sigma_c$ we have

$$(1 - a_2 e^{-\lambda_2 s}) \zeta_{A,\Lambda}(s) = 1 + a_3 e^{-\lambda_3 s} + a_5 e^{-\lambda_5 s} + \ldots$$

where in the right hand side all the terms with even subscript have cancelled since $a_2 a_k = a_{2k}$ and $\lambda_2 + \lambda_k = \lambda_{2k}$.

Similarly, if we multiply term by term this equality by $(1 - a_3 e^{-\lambda_3 s})$, all the terms with subscript $3k$ from the right hand side will cancel and so on. After a finite number of steps in which we use consecutive prime subscripts, given an arbitrary $n$, the right hand side will contain only terms with the subscripts greater than $n$. It can be easily checked that $\lim_{n \to \infty} \Sigma'_{k \geq n} a_k e^{-\lambda_k s} = 0$ uniformly on compact sets of $\operatorname{Re} s > \sigma_c$.

Here $\Sigma'$ means that we use only subscripts $k$ which are not multiples of $p$, $p < n$. Indeed, there is a constant C such that

$$|\Sigma'_{k \geq n} a_k e^{-\lambda_k s}| \leq \Sigma_{k \geq n} |a_k| e^{-\lambda_k (\sigma - \sigma_0)} e^{-\lambda_k \sigma_0} =$$

$$e^{-\lambda_n (\sigma - \sigma_0)} \Sigma_{k \geq n} |a_k| e^{-(\lambda_k - \lambda_n)(\sigma - \sigma_0)} e^{-\lambda_k \sigma_0} \leq C e^{-\lambda_n \sigma},$$

which tends to zero uniformly on compact sets of $\operatorname{Re} s > \sigma_c$. The last inequality takes place since $e^{-(\lambda_k - \lambda_n)(\sigma - \sigma_0)} \leq 1$, $\Sigma_{k \geq n} |a_k| e^{-\lambda_k \sigma_0}$ is convergent and $e^{-\lambda_n \sigma_0}$ is a constant independent of $k$. The conclusion is that (5) holds for $\operatorname{Re} s > \sigma_c$.

## 3. FUNCTIONS SATISFYING A RIEMANN-TYPE OF FUNCTIONAL EQUATION BUT NOT SATISFYING RH

The set of functions obtained by analytic continuation of general Dirichlet series and which satisfy a Riemann-type of functional equation can be partitioned into two subsets: one including all those functions which satisfy RH and the other one including all the functions of this type which do not satisfy RH. Since the definition of the non trivial zeros requires that the function satisfies a Riemann-type of functional equation, only functions belonging to this class make sense to be taken into consideration. We need to emphasize this fact, since sometimes Hurwitz Zeta functions, which do not satisfy a Riemann-type of functional equation, or some elementary functions satisfying such an equation, but not being Dirichlet series, are brought into context. The off critical line non trivial zeros found for the Davenport and Heilbronn function show that the second class is not empty. We have shown in [11] that infinitely many functions of that type can be obtained by a similar construction. It is then interesting to put into evidence some general properties of these functions. The following theorem serves this purpose.

**Theorem 2**. *Suppose that $f(s)$ is obtained by analytic continuation to the whole complex plane of a general Dirichlet series and it satisfies a Riemann-type of functional equation. If $f(s)$ does not satisfy RH, then:*

(a) *For every two distinct non trivial zeros $s_1 = \sigma + it$ and $s_2 = 1 - \sigma + it$ there is $\tau_0$, $0 < \tau_0 < 1$ such that $f'(s(\tau_0)) = 0$, where $s(\tau) = (1-\tau)s_1 + \tau s_2$, i.e. the derivative cancels at a point belonging to the segment $I$ determined by $s_1$ and $s_2$.*

(b) $\operatorname{Re} f(s(\tau_0)) < 1/2$

*Proof:* Let $z = z(\tau) = f(s(\tau))$, $0 \leq \tau \leq 1$ be the image of the segment $I$ by the function $f(s)$. It is a closed curve $\gamma$ passing through the origin since $f(s(\tau)) = 0$ when $\tau = 0$ and when $\tau = 1$. As $f(s)$ is analytic and $s(\tau)$ is differentiable, we can differentiate with respect to $\tau$ for $0 < \tau < 1$ and we get $z'(\tau) = f'(s(\tau))s'(\tau) = f'(s(\tau))(s_2 - s_1) = (1 - 2\sigma)f'(s(\tau))$. If $f'(s)$ does not cancel at any point of $I$, then $z'(\tau) \neq 0$, $0 < \tau < 1$, hence the curve $\gamma$ has a tangent at every one of its points, except the origin. We can suppose that $1 - 2\sigma > 0$, thus $\operatorname{Arg} z'(\tau) = \operatorname{Arg} f'(s(\tau))$, i.e. the tangent to $\gamma$ at $z(\tau)$ has the same orientation as the position vector of $f'(s(\tau))$. Let $s_0$ be the zero of $f'(s)$ determined by $s_1$ and $s_2$ in the sense that components containing $s_1$ and respective $s_2$ of the pre-image of a circle centered at the origin touch each other at $s_0$. Let us denote by $L$ the ray from the origin passing through $f(s_0)$. Continuation by $f(s)$ over $L$ starting from $s_1$ and $s_2$ gives rise to two unbounded curves $\eta_1$ and $\eta_2$ passing through $s_0$. The unboundedness of $\eta_1$ and of $\eta_2$ is guaranteed by the fact that the continuations are unlimited, since there is no singular point in their way. If we ignore the parts of $\eta_1$ and $\eta_2$ from $s_1$ and $s_2$ to $s_0$, what remains from $\eta_1 + \eta_2$ is an unbounded curve $\eta$ separating $s_1$ and $s_2$. Then $I$ must intersect $\eta$, which means that $\gamma$ must intersect the ray $L$. But then $I$ must intersect both $\eta_1$ and $\eta_2$, which is possible only if $I$ passes through $s_0$. This proves the part (a) of the theorem and shows at the same time that $\gamma$ passes through $f(s_0)$. Let $\tau_0$ be such that $s(\tau_0) = s_0$. Then, since $f'(s(\tau_0)) = 0$, we conclude that $z'(\tau_0) = 0$, which means that $\gamma$ has at $f(s_0)$ a turning back point.

Let $\Omega$ and $\Omega'$ be fundamental domains of $\zeta_{A,\Lambda}(s)$ such that $\sigma + it \in \Omega$, ($\sigma \geq 1/2$) and $1 - \sigma + it \in \Omega'$. Looking at Fig. 2 might be helpful, although any other configuration of fundamental domains containing $\sigma + it$ and $1 - \sigma + it$ is acceptable. The function $\zeta_{A,\Lambda}(s)$ maps conformally $\Omega$ onto the complex plane with a slit $L$ and it maps conformally $\Omega'$ onto the complex plane with a slit $L'$ such that $0 \notin L \cup L'$. Let $B \subseteq \Omega$ and $B' \subseteq \Omega'$ be components of the pre-image by $\zeta_{A,\Lambda}(s)$ of $\mathbb{C} \setminus (L \cup L')$. Then $\sigma + it \in B$, $1 - \sigma + it \in B'$. Let us suppose that $\Omega$ is the domain unbounded to the right (one of $\Omega$ and $\Omega'$ has always this property) and define $\phi : B \rightarrow B'$ by the formula

$$\phi(s) = \zeta_{A,\Lambda}^{-1}|_{B'} \circ \zeta_{A,\Lambda}|_B (s) \tag{6}$$

The function $\phi(s)$ is correctly defined, since $\zeta_{A,\Lambda}(s)$ is conformal, hence injective in every fundamental domain. The function $\phi(s)$ can be continued analytically to an involution on the simply connected domain $B \cup \overline{B'}$. The equality $\phi(\phi(s)) = s$ implies $\phi'(\phi(s))\phi'(s) \equiv 1$, thus $|\phi'(s)| < 1$ if and only if $|\phi'(\phi(s))| > 1$.

On the other hand, $\zeta_{A,\Lambda}(\phi(s)) = \zeta_{A,\Lambda}(s)$ implies that $\zeta'_{A,\Lambda}(\phi(s))\phi'(s) = \zeta'_{A,\Lambda}(s)$, thus $|\zeta'(\phi(s))| < |\zeta'_{A,\Lambda}(s)|$ if and only if $|\phi'(s)| < 1$, in other words $\zeta_{A,\Lambda}$ performs locally a contraction in $B$ iff it performs an expansion in $B'$. Since the pre-image by $\zeta_{A,\Lambda}(s)$ of the interval $(0,1)$ has a finite component in $B'$ and an infinite one in $B$ we deduce that $\zeta_{A,\Lambda}(s)$ is a contraction in $B$.

Let us denote by $h_2$ and $h_1$ the corresponding inverse functions and notice that $h_2$ is a local expansion, while $h_1$ is a local contraction. Let $C$ be the arc mapped by $h_2$ onto the segment $\delta_2$ determined by $s_0$ and $s_2$. The arc $C$ is also the image by $\zeta_{A,\Lambda}(s)$ of an arc $\delta_1$ (not necessarily a segment) connecting $s_1$ and $s_0$. The length of the segment $\delta_2$ is $\int_C |h'_2(z)||dz|$, while the length of the arc $\delta$ is $\int_C |h'_1(z)||dz|$ and the first is greater than the second, in particular it is greater than the length of the segment between $s_1$ and $s_0$, which proves the part (b) of the theorem.

## 4. THE PROOF OF THE RH FOR A CLASS OF FUNCTIONS

One would expect this section to be long and to contain very intricate arguments. In fact, it is short and elementary and probably proves more than GRH.

**Theorem 3**. *Suppose that $n \rightarrow a_n$ is a totally multiplicative function and $n \rightarrow \lambda_n$ satisfies the condition* (4). *If $\sigma_c \leq 1/2$ and $\zeta_{A,\Lambda}(s)$ satisfies the equation* (3), *then for every non trivial zero $\sigma + it$ of $\zeta_{A,\Lambda}(s)$ we have $\sigma = 1/2$.*

*Proof:* We notice first that the hypothesis of this theorem is fulfilled by the Dirichlet L-functions for which $\chi(n)$ is a totally multiplicative function, $\lambda_n = \ln n$, and $\sigma_c = 0$. Then the results from [12] and [14] are valid, as long as the zero $\sigma + it$ does not belong to a curve $\Gamma_{k,0}$. Hence we could deal in this theorem only with the case where $\sigma + it \in \Gamma_{k,0}$. One can follow the reasoning looking at Fig. 2, although everything we are saying is true for any couple of non trivial zeros symmetric with respect to the critical line.

The function $\phi(s)$ defined by (6) is such that $\phi'(s) \neq 0$ in $B \cup \overline{B'}$ (there are no branch points of $\phi(s)$ in $B \cup \overline{B'}$, where the univalence would be violated) and $\phi(\sigma + it) = 1 - \sigma + it$.
By the definition of $\phi(s)$, we have that $\zeta_{A,\Lambda}(\phi(s)) = \zeta_{A,\Lambda}(s)$. Let us denote $\Omega_c = (B \cup B') \setminus [\{s \mid \mathrm{Re}\, s \leq \sigma_c\} \cap \{s \mid \mathrm{Re}\, \phi(s) \leq \sigma_c\}]$ and define $f_n(s)$ in $\Omega_c$ as follows:
$$f_n(s) = \Pi_{p \leq n}(1 - a_p e^{-\lambda_p s}), \quad in \ (B \cup B') \cap \{s \mid \mathrm{Re}\, s > \sigma_c\} \qquad (7)$$

$$f_n(s) = \Pi_{p \leq n}(1 - a_p e^{-\lambda_p \phi(s)}), \quad in \ (B \cup B') \cap \{s \mid \mathrm{Re}\, s \leq \sigma_c\} \cap \{s \mid \mathrm{Re}\, \phi(s) > \sigma_c\}$$

where $p$ runs through prime numbers. The functions $f_n(s)$ are analytic in $\Omega_c$, which includes neighborhoods of $\sigma + it$ and $1 - \sigma + it$. By (5) and the fact that $\zeta_{A,\Lambda}(\phi(s)) = \zeta_{A,\Lambda}(s)$, we have

$$\lim_{n \to \infty}[1 / f_n(s)] = \zeta_{A,\Lambda}(s) \qquad (8)$$

uniformly on compact sets included in $\Omega_c$. Let us notice that by (2), the abscissa of convergence $\sigma_c$ depends only on $a_p$ and $\lambda_p$ and it is such that for $\mathrm{Re}\, s > \sigma_c$, due to (7) and (8), the equation $1 - a_p e^{-\lambda_p s} = 0$ cannot have any solution $s_0 \in \Omega_c$ for any prime number $p$, since otherwise there would be $n_0$ such that for $n \geq n_0$, $f_n(s_0) = 0$ hence $\infty = \lim_{n \to \infty} [1/f_n(s_0)] = \zeta_{A,\Lambda}(s_0)$, which is absurd. Then, for $s \in \Omega_c$, $s \neq \sigma + it$ we have

$$\lim_{n \to \infty} f_n(\phi(s)) / f_n(s) = \zeta_{A,\Lambda}(s) / \zeta_{A,\Lambda}(\phi(s)) = 1 \qquad (9)$$

uniformly on compact subsets of $\Omega_c \setminus \{\sigma + it\}$. Yet, the functions $f_n(\phi(s)) / f_n(s)$ are defined throughout $\Omega_c$ and they are analytic in $\Omega_c$. Indeed, the denominator cannot cancel in $\Omega_c$, as we have seen. By the functional equation, $\sigma + it$ and $1 - \sigma + it$ are zeros of $\zeta_{A,\Lambda}(s)$ of the same order of multiplicity, hence the function $\Phi(s) = \zeta_{A,\Lambda}(s) / \zeta_{A,\Lambda}(\phi(s))$, has a removable singularity at $\sigma + it$, i.e. $\Phi(s) \equiv 1$ in $\Omega_c$. Due to (8), we can say that $f_n(\phi(s)) / f_n(s)$ converges uniformly to 1 on compact subsets of $\Omega_c$

.

For $s = \sigma + it$ the limit of this sequence is:

$$\lim_{n\to\infty} \Pi_{p\leq n}[(1 - a_p e^{-\lambda_p(1-\sigma+it)}) / (1 - a_p e^{-\lambda_p(\sigma+it)})] =$$

$$\lim_{n\to\infty} \Pi_{p\leq n}[(e^{\lambda_p(1-\sigma+it)} - a_p) / (e^{\lambda_p(\sigma+it)} - a_p)]e^{\lambda_p(2\sigma-1)} \quad (10)$$

If $2\sigma - 1 > 0$ we have that $e^{\lambda_p(2\sigma-1)} \to +\infty$ as $p \to \infty$. The condition $\sigma_c \leq 1/2$ was necessary in order to make sure that the infinite product (5) and therefore the sequence (9) converges in a neighborhood of $\sigma + it$, where $\sigma > 1/2$. Yet, this last condition implies that the limit (10) cannot exist, since the numerator in (10) cannot cancel in order to offset the factor $e^{\lambda_p(2\sigma-1)}$. This contradicts the established fact that the sequence (9) converges uniformly to 1 on compact subsets of $\Omega_c$. Therefore $\sigma = 1/2$, in which case the limit (10) is 1. Thus, for every non trivial zero $\sigma + it$ of $\zeta_{A,\Lambda}(s)$ we have $\sigma = 1/2$, i.e. the RH is true for this class of functions.

We notice again that this class of functions includes the Dirichlet L-functions, therefore Theorem 3 completes the proof of Theorem 9 from [12], where a different geometric method has been used in which the possibility of the violation of RH due to a zero on $\Gamma_{k,0}$ has been overlooked.